\documentclass[letterpaper, 10 pt, conference]{ieeeconf}  

\IEEEoverridecommandlockouts                              

\usepackage{graphics} 
\usepackage{epsfig} 
\usepackage{times} 
\usepackage{amsmath} 
\usepackage{amssymb}  
\usepackage{graphicx}      
\usepackage{booktabs}
\usepackage{mathtools}

\title{\LARGE \bf
A Bayesian Neural ODE for a Lettuce Greenhouse
}
\author{Sjoerd Boersma and Xiaodong Cheng
\thanks{Sjoerd Boersma and Xiaodong Cheng are both with the Mathematical and Statistical Methods Group (Biometris), Department of Plant Science, Wageningen University \& Research,
			6700 AA Wageningen, The Netherlands.
        {\tt\small sjoerd.boersma@wur.nl, xiaoding.cheng@wur.nl}}
}

\begin{document}

\maketitle
\thispagestyle{empty}
\pagestyle{empty}

\begin{abstract}
Greenhouse production systems play a crucial role in modern agriculture, enabling year-round cultivation of crops by providing a controlled environment. However, effectively quantifying uncertainty in modeling greenhouse systems remains a challenging task. In this paper, we apply a novel approach based on sparse Bayesian deep learning for the system identification of lettuce greenhouse models. The method leverages the power of deep neural networks while incorporating Bayesian inference to quantify the uncertainty in the weights of a Neural ODE.  The simulation results show that the generated model can capture the intrinsic nonlinear behavior of the greenhouse system with probabilistic estimates of environmental variables and lettuce growth within the greenhouse.

\end{abstract}
\section{INTRODUCTION}

Greenhouses have become one of the most rapidly expanding  production systems worldwide in the agricultural industry. They offer a controlled or partially closed environment for plants, protecting them against external weather conditions while providing essential water and nutrients. For better environmental and economic sustainability, effective management of modern greenhouses needs to balance multiple factors, such as temperature, humidity, light density, and nutrient levels, which are regulated by indoor climate control actuators, for example, the heating, ventilation, and carbon dioxide supply. From the perspective of control design, the operation of these actuators usually leans on accurate and reliable mathematical models, allowing the prediction of both greenhouse climate and crop growth \cite{katzin2022process}. Nevertheless, due to sensor noises and environmental fluctuations, system uncertainties are inevitable in greenhouse modeling and control, and it may lead to inaccurate predictions and erroneous control actions if the system uncertainties are not appropriately considered \cite{Henten2003,Svensen2024cc}. 
  
The conventional way of greenhouse modeling is often deterministic and based on first principles, which means that the dynamic greenhouse climate and crop growth are modeled using the scientific knowledge of the physical and biological processes \cite{bakker1995greenhouse,vanHenten1991,vanHenten1994}.
This method tends to produce complex models with intricate nonlinear dynamics and a large number of model parameters. Therefore, the first-principle models usually lack of flexibility to account for different greenhouse settings or scenarios. Moreover, the first principle approaches are usually difficult to effectively capture the uncertainty of the system, as they may demand expensive computation and substantial amounts of data to statistically determine the uncertainty of the model parameters. 

On the other hand, with the recent advances in machine learning, different data-driven techniques have been employed to discover models from the data that are collected online from the greenhouse sensors and actuators, see e.g., \cite{wang2009svm,dariouchy2009prediction,taki2018applied,Chen2022}, where the support vector machine (SVM) and neural networks are applied to model the indoor greenhouse environment. The data-driven methods, compared to the first-principle approaches, are more flexible for adapting to different greenhouse settings, as they can learn from the data without relying on the prior knowledge of the system. However, this method also has some challenges, such as the need for large amounts of data and computational power to obtain reliable and accurate models. For instance, the data-driven models may not be able to generalize well to new situations or cope with the noise or errors in the data.

The recently emerged method, called Bayesian deep learning, appears to be a promising approach to address the above challenges in modeling greenhouse dynamics. In contrast to the traditional deep learning models, such as neural networks, that typically rely on deterministic architectures and point estimates for parameters, Bayesian deep learning appears as a probabilistic modeling approach accounting for uncertainty in the model parameters. More preciously, instead of providing a single fixed value for each parameter, Bayesian deep learning models give a probability distribution over possible values for each parameter. We refer to \cite{wang2020survey} for an overview. The Bayesian framework offers a probabilistic perspective in interpreting systems so that model parameters and prediction uncertainties can be quantified with probabilities. Recently, this methodology has also been combined with neural ordinary differential equations (ODEs), to learn the dynamics of a physical system with the quantification of
the uncertainty in the weights of a Neural ODE \cite{dandekar2020bayesian}. 

In this paper, we apply the Bayesian deep learning approach to modeling greenhouses. We take a lettuce greenhouse discussed in \cite{vanHenten1994} as a benchmark and aim to build a Bayesian neural ODE model that can predict the environmental variables and lettuce growth in a greenhouse with a measure of the model uncertainty. The neural network structure is the Multi-Layer Perceptron (MLP) and a sparse regression method described in \cite{zhou2022sparse} is used to enforce a relatively sparse structure in the neural network. By treating model parameters as probability distributions, we can quantify uncertainty and make informed predictions of the environmental variables and crop growth in a greenhouse. 
\clearpage
{\it Outline}: The paper is organized as follows: Section~\ref{sec:model} describes the setup of a lettuce greenhouse and the data generation process, and Section~\ref{sec:neuralOde} presents the main method of the Bayesian neural ODE and relevant the optimization problem.  Section~\ref{sec:results} shows the simulation results to validate the generated model, and Section~\ref{sec:conclusions} provides conclusion remarks of the paper.

\section{LETTUCE GREENHOUSE MODEL}
\label{sec:model}

This work exploits a benchmark model for generating training and validation data. This nonlinear lettuce greenhouse model has been published in~\cite{vanHenten1994} and can be written in state-space form as follows:
\begin{equation}
	\begin{aligned}
		x(k+1) &= f\big(x(k),u(k),d(k),p\big), \\
		y(k)   &= x(k), 
		\label{eq:system}
	\end{aligned}
\end{equation}
with discrete time $k \in \mathbb{Z}^{0+}$, state $x(k),y(k) \in \mathbb{R}^4$, controllable input $u(k) \in \mathbb{R}^3$, weather disturbance $d(k) \in \mathbb{R}^4$, parameter $p = (p_1~p_2~\dots~p_{28})^T \in \mathbb{R}^{28}$ and $f(\cdot)$ a nonlinear function that is, together with initial conditions, given in the Appendix. Table~\ref{tab:signals} provides the meaning of the signals.

\begin {table}[h]
\caption{Meaning of the state $x(t)$, control signal $u(t)$ and disturbance $d(t)$.}
\begin{center}
	\begin{tabular}{l  l }
		\toprule 			
		$x_{1}(t)$ & dry-weight (kg/m$^2$)  \\
		$x_{2}(t)$ & indoor CO$_2$ (kg/m$^3$)  \\
		$x_{3}(t)$ & indoor temperature (deg $^\circ$) \\
		$x_{4}(t)$ & indoor humidity (kg/m$^3$)  \\	    
		$d_{1}(t)$ &  radiation (W/m$^2$) \\
		$d_{2}(t)$ &  outdoor CO$_2$ (kg/m$^3$) \\
		$d_{3}(t)$ & outdoor temperature (deg $^\circ$) \\
		$d_{4}(t)$ & outdoor humidity (kg/m$^3$)\\
		$u_{1}(t)$ &  CO$_2$ injection (mg/m$^2$/s) \\
		$u_{2}(t)$ &  ventilation (mm/s) \\
		$u_{3}(t)$ &  heating (W/m$^2$) \\
		\bottomrule
	\end{tabular}
	\label{tab:signals}
\end{center}
\end {table}

\begin{figure}[!h]
	\centering
	\includegraphics[width=.49\textwidth]{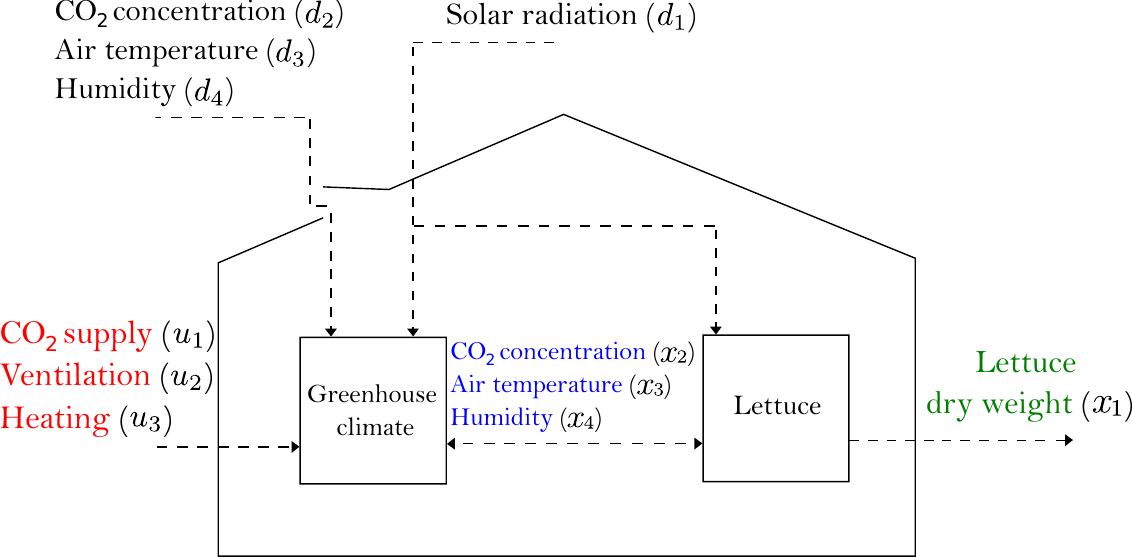}
	\caption{Schematic representation of lettuce greenhouse model. The arrows indicate interactions between the greenhouse-lettuce model with the control signal $u(k)$, weather disturbance $d(k)$ and measurement $x(k)$. \label{fig:greenhouse}}
\end{figure}

Figure~\ref{fig:greenhouse} provides a graphical representation of the data generating lettuce greenhouse model.

The model in~\eqref{eq:system} is used as the simulation model to generate data that is used for training the Bayesian neural ODE. The data generating procedure is described in the following section.

\subsection{Data Generation}

Providing different weather disturbances and control inputs, $N_s$ scenarios are generated $N_s$ scenarios are generated in total, each spanning $N$ time steps. The generated data in each scenario at time step $k$ is defined as a throuple:
\begin{equation}
	X^j(k) = \{ x^j(k),u^j(k),d^j(k) \},
\end{equation} 
for $\quad j=1,\ldots,N_s$ and $k=1,\ldots,N$, where $d^j(k)$ comes from real weather data as described in~\cite{Kempkes2014} and is collected at a greenhouse in Holland. The control inputs are designed depending on the weather conditions: 
\begin{equation}
	\begin{aligned} 
		u_1^j(k)    &= d_1^j(k)/10, \\
		u_2^j(k)    &= \begin{cases} \mbox{0,} & \mbox{if } d_1^j(k) > 1 \\ \mbox{0.1,} & \mbox{otherwise} \end{cases}, \\
		u_3^j(k)    &= 20+d_1^j(k)/5.
	\end{aligned} 
\end{equation}
That is, both the CO$_2$ injection and heating are controlled proportionally to the radiation levels. Consequently, additional CO$_2$ is injected and the heating is cranked up during the day. This relatively straightforward control law uses the fact that lettuce growth is primarily active during daylight hours and additional CO$_2$ and heating can thereby actually contribute to the growth. There is no ventilation when the sun is shining since then, photosynthesis takes place, and additional CO$_2$ must not be lost via the ventilation. However, fresh air enters the greenhouse when there is almost no radiation since then, $u_2^j(k)=0.1$.  

\begin{figure*}
	\centering
	\includegraphics[width=\textwidth]{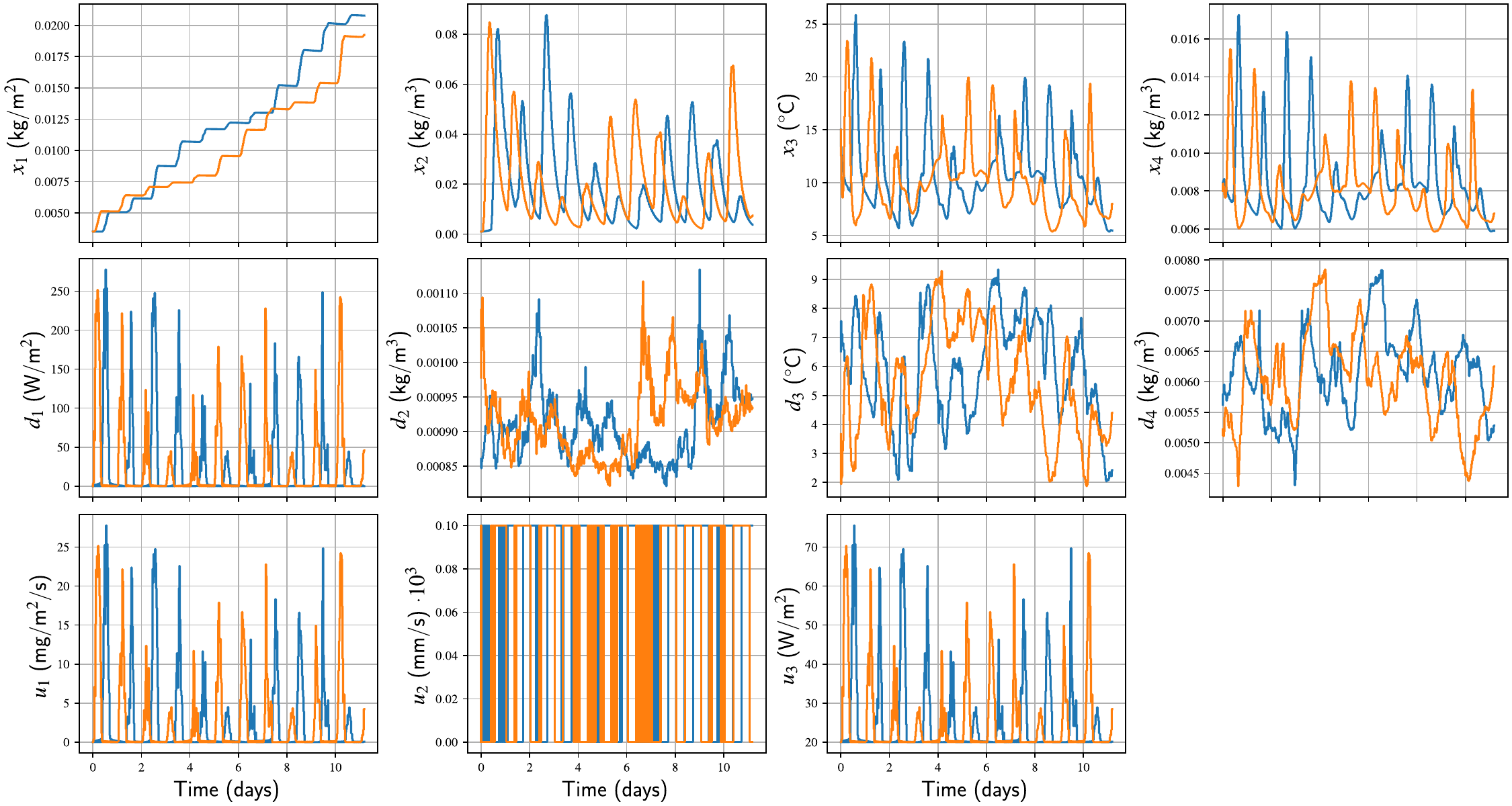}
	\caption{Two example scenarios $X^1(k),X^2(k)$ that are, after normalization, used in the training data set. In total, $N_s$ of these scenarios are used for training the Bayesian neural ODE. \label{fig:scenarios}}
\end{figure*}

The variation in the inputs imposes variation in the data since $x^j(k)$ will be different for each scenario. Figure~\ref{fig:scenarios} depicts two example scenarios. In the following section, the data from the multiple scenarios will be used to learn a Bayesian Neural ODE that captures the dynamics of the lettuce greenhouse system.


\section{METHODOLOGY}
\label{sec:neuralOde}

The used methodology is taken from~\cite{zhou2022sparse}. Here, a sparse Bayesian deep learning method for dynamic system identification is presented, which allows identifying (training) a dynamical model with uncertain parameters. Such a stochastic model is suitable for the application at hand due to the uncertain nature of the weather, which influences the evolution of state variables. So in order to predict the future state variables, the uncertain future weather needs to be accounted for and one way to do this is via the use of a stochastic model as done in this work.

\subsection{Bayesian Neural ODE}

The employed method is used in this work to identify a model having the following structure:

\begin{equation}
	\begin{aligned}
		\bar{\dot{\hat{x}}}(t) &= \hat{f}\big(\bar{\hat{x}}(t),\bar{u}(t),\bar{d}(t),\hat{p}\big), \\
		\bar{\hat{y}}(t)       &= \bar{\hat{x}}(t), 
		\label{eq:model}
	\end{aligned}
\end{equation}
with $\bar{\cdot}$ indicating normalized variables (explained later) and continuous time $t \in \mathbb{R}$. The vector of the model parameters $\hat{p}$ belongs to a multivariate normal distribution with mean $\mu_{\hat{p}}$ and standard deviation $\sigma_{\hat{p}}$, i.e.,
\begin{equation}
	\hat{p} \sim \mathcal{N}(\mu_{\hat{p}},\sigma_{\hat{p}}).
\end{equation}
The individual parameters itself are assumed to be uncorrelated so $\sigma_{\hat{p}}=\texttt{diag}(\sigma_{\hat{p}_1},\ldots,\sigma_{\hat{p}_{n_{\hat{p}}}})$. The nonlinear function $\hat{f}(\cdot)$ is defined by layers, activation functions, and neurons. The input and output layers are linear. A user-defined number of hidden layers can be placed in between the input and output layer with each its own number of neurons and activation functions. These (number of hidden layers, number of neurons, and activation functions) are tuning variables. For a more detailed description of the structure $\hat{f}(\cdot)$, the reader is referred to~\cite{zhou2022sparse}.

\subsection{Feature and Target Definitions}

In~\eqref{eq:approxDer}, the derivative is defined as the target variable. Since this is generally not available in practice, the approximate derivative is defined as the target variable for each scenario $j$, i.e.,
\begin{equation}
	\dot{x}^j(k) \coloneqq \frac{x^j(k+1)-x^j(k)}{h},
	\label{eq:approxDer}
\end{equation} 
with sample period $h$. Consequently, the feature variables are $x^j(k),u^j(k),d^j(k)$. Both the target and features are normalized as follows:
\begin{equation}
\begin{aligned}
    \bar{\dot{x}}^j(k) &= (\dot{x}^j(k)-\mu_{\dot{x}})./\sigma_{\dot{x}}, \\
    \bar{x}^j(k) &= (x^j(k)-\mu_{x})./\sigma_{x}, \\
    \bar{d}^j(k) &= (d^j(k)-\mu_{d})./\sigma_{d}, \\
    \bar{u}^j(k) &= (u^j(k)-\mu_{u})./\sigma_{u},    
\label{eq:normalization}
\end{aligned}
\end{equation} 
with $./$ element wise division, $\mu_{\dot{x}} \in \mathbb{R}^4, \mu_{x} \in \mathbb{R}^4,\mu_{d} \in \mathbb{R}^4,\mu_{u} \in \mathbb{R}^3$ the sample means of $\dot{x}^j(k), x^j(k), d^j(k), u^j(k)$ over all $j$ and $k$, respectively. Furthermore, $\sigma_{\dot{x}} \in \mathbb{R}^4,\sigma_{x} \in \mathbb{R}^4, \sigma_{d} \in \mathbb{R}^4, \sigma_{u} \in \mathbb{R}^3$ are the sample standard deviations of the time series data $x(k), d(k), u(k)$, respectively.   

These definitions yield the following data matrices that are used for training the Bayesian neural ODE:
\begin{equation} \scriptsize
	\begin{aligned}
		X = \underbrace{\begin{pmatrix} 
							  \bar{x}^1(0)       & \bar{u}^1(0)       & \bar{d}^1(0) \\
			                  \bar{x}^1(1)       & \bar{u}^1(1)       & \bar{d}^1(1) \\
							  \vdots       & \vdots       & \vdots \\			
							  \bar{x}^{1}(N)     & \bar{u}^1(N)     & \bar{d}^1(N) \\ 
						      \bar{x}^2(0)       & \bar{u}^2(0)       & \bar{d}^2(0) \\
							  \vdots       & \vdots       & \vdots \\	
						      \bar{x}^{N_s}(0)   & \bar{u}^{N_s}(0)   & \bar{d}^{N_s}(0) \\
						      \vdots       & \vdots       & \vdots \\
						      \bar{x}^{N_s}(N)   & \bar{u}^{N_s}(N) & \bar{d}^{N_s}(N)
					        \end{pmatrix}}_{\texttt{features}}, \quad 
				        		Y = \underbrace{\begin{pmatrix} 
				        		\bar{\dot{x}}^1(0)        \\
				        		\bar{\dot{x}}^1(1)        \\
				        		\vdots        \\			
				        		\bar{\dot{x}}^{1}(N)    \\ 
				        		\bar{\dot{x}}^2(0)        \\
				        		\vdots        \\	
				        		\bar{\dot{x}}^{N_s}(0)    \\
				        		\vdots        \\
				        		\bar{\dot{x}}^{N_s}(N) 
				        \end{pmatrix}}_{\texttt{target}}
	\end{aligned}
 \label{eq:data}
\end{equation}

The state derivative~\eqref{eq:approxDer} is used as the target because it arguably contains more information than the state itself in unstable systems. Indeed, the dry-weight $x_1(k)$ is an unstable state variable (see Fig.~\ref{fig:scenarios}). 

\subsection{Optimization Problem}

The training of the neural network boils down to the following optimization problem:
\begin{equation}
	\min_{\mu_{\hat{p}},\sigma_{\hat{p}}} \quad \frac{1}{N_s N} \sum_{j=0}^{N_s} \sum_{k=0}^{N} \big( \bar{\dot{x}}^j(k) - \bar{\hat{\dot{x}}}(k) \big)^T \big( \bar{\dot{x}}^j(k) - \bar{\hat{\dot{x}}}(k) \big),
	\label{eq:optimizationNN}
\end{equation} 
with $N$ the number of used data points in the identification procedure, $N_s$ the number of scenarios taken into account by the identification, $\bar{\hat{\dot{x}}}(k)$ the normalized state generated by~\eqref{eq:model} hence depends on the decision variable (model parameters) $\hat{p} \sim \mathcal{N}(\mu_{\hat{p}},\sigma_{\hat{p}})$ and $\bar{\dot{x}}^j(k)$ is the normalized state coming from~\eqref{eq:system}, the true system. The optimization problem defined in~\eqref{eq:optimizationNN} is solved using the code provided by~\cite{zhou2022sparse}.

\subsection{Simulating the Bayesian Neural ODE}

In order to simulate the stochastic model in~\eqref{eq:model}, it is first temporally discretized. Euler's method is used to do this:
\begin{equation}
	\begin{aligned}
		\bar{\hat{x}}(k+1)     &= \bar{\hat{x}}(k) + h \cdot \hat{f}\big(\bar{\hat{x}}(k),\bar{u}(k),\bar{d}(k),\hat{p}\big), \\
		\bar{\hat{y}}(k)       &= \bar{\hat{x}}(k), \qquad \hat{p} \sim \mathcal{N}(\mu_{\hat{p}},\sigma_{\hat{p}}), 
		\label{eq:modelDiscrete}
	\end{aligned}
\end{equation}
where the approximation $\bar{\dot{\hat{x}}}(t) \approx \big( \bar{\hat{x}}(k+1)- \bar{\hat{x}}(k) \big)/h$ is exploited to go from~\eqref{eq:model} to~\eqref{eq:modelDiscrete}. Since $\hat{p}$ is stochastic and $\hat{f}(\cdot)$ is nonlinear, $\hat{p}$ needs to be sampled yielding the following model description:
\begin{equation}
	\begin{aligned}
		\bar{\hat{x}}^i(k+1)     &= \bar{\hat{x}}^i(k) + h \cdot \hat{f}\big(\bar{\hat{x}}^i(k),\bar{u}(k),\bar{d}(k),\hat{p}^i\big), \\
		\bar{\hat{y}}^i(k)       &= \bar{\hat{x}}^i(k), \qquad \text{for}~i=1,\ldots,N_{\hat{p}},
		\label{eq:modelDiscrete1}
	\end{aligned}
\end{equation}
where $\hat{p}^i$ is a deterministic parameter vector drawn from the normal distribution with mean $\mu_{\hat{p}}$ and standard deviation $\sigma_{\hat{p}}$. So the deterministic model in~\eqref{eq:modelDiscrete1} is simulated $N_{\hat{p}}$ times. Then, $\hat{x}^i(k),u(k),d(k)$ is obtained using the inverse of~\eqref{eq:normalization}, i.e. the model's estimation is denormalized. Consequently, a mean state and confidence interval are computed and plotted.

\section{SIMULATION RESULTS}
\label{sec:results}

\subsection{Weather Data}
The assumed to be measured weather data $d(k)$ used in this work is real data and described in~\cite{Kempkes2014}. The origin of the data is ``the Venlow Energy greenhouse'', which is located in Bleiswijk, Holland. The collected data points are re-sampled to a sample period of $h=30$ minutes.

\subsection{Simulation Setup}
The identification method is tested following the steps:
\begin{enumerate}
    \item generate $N_s$ scenarios (data) of each 11 days long using~\eqref{eq:system},
    \item normalize the data and create the data matrices with target and features~\eqref{eq:data},
    \item solve the optimization problem in~\eqref{eq:optimizationNN},
    \item simulate $N_{\hat{p}}$ times the identified neural ODE in~\eqref{eq:model} from day 11 until 14 using~\eqref{eq:modelDiscrete1},
    \item de-normalize the state variable obtained in the previous step.
    \item compute the mean state variable from the previous step and its confidence intervals.
\end{enumerate}
Remember that in the first step, each generated scenario is different since different weather data is used as an input for each scenario. Also note that the training data is from day 0 to day 11, while the validation data is from day 11 to day 14. So the neural network predicts the upcoming 4 days given the past 11 days. In practice, this can be applied in a recurrent way and each day a newly trained neural ODE can be trained to predict the upcoming 3 days (or potentially shorter/longer).  

\subsection*{Neural Network Identification Settings}
The Bayesian neural network has one linear input and one linear output layer. In addition to that, the network has one hidden layer with four neurons. These selected numbers of hidden layers and neurons are chosen such that an acceptable validated model can identified that strikes a balance between overfitting and predictive capability. While reducing the number of layers and neurons prevents overfitting, using fewer neurons will significantly reduce the predictive capabilities of the deep neural network.

\subsection{Identification Results}

\begin{figure}[!h]
	\centering
	\includegraphics[width=.47\textwidth]{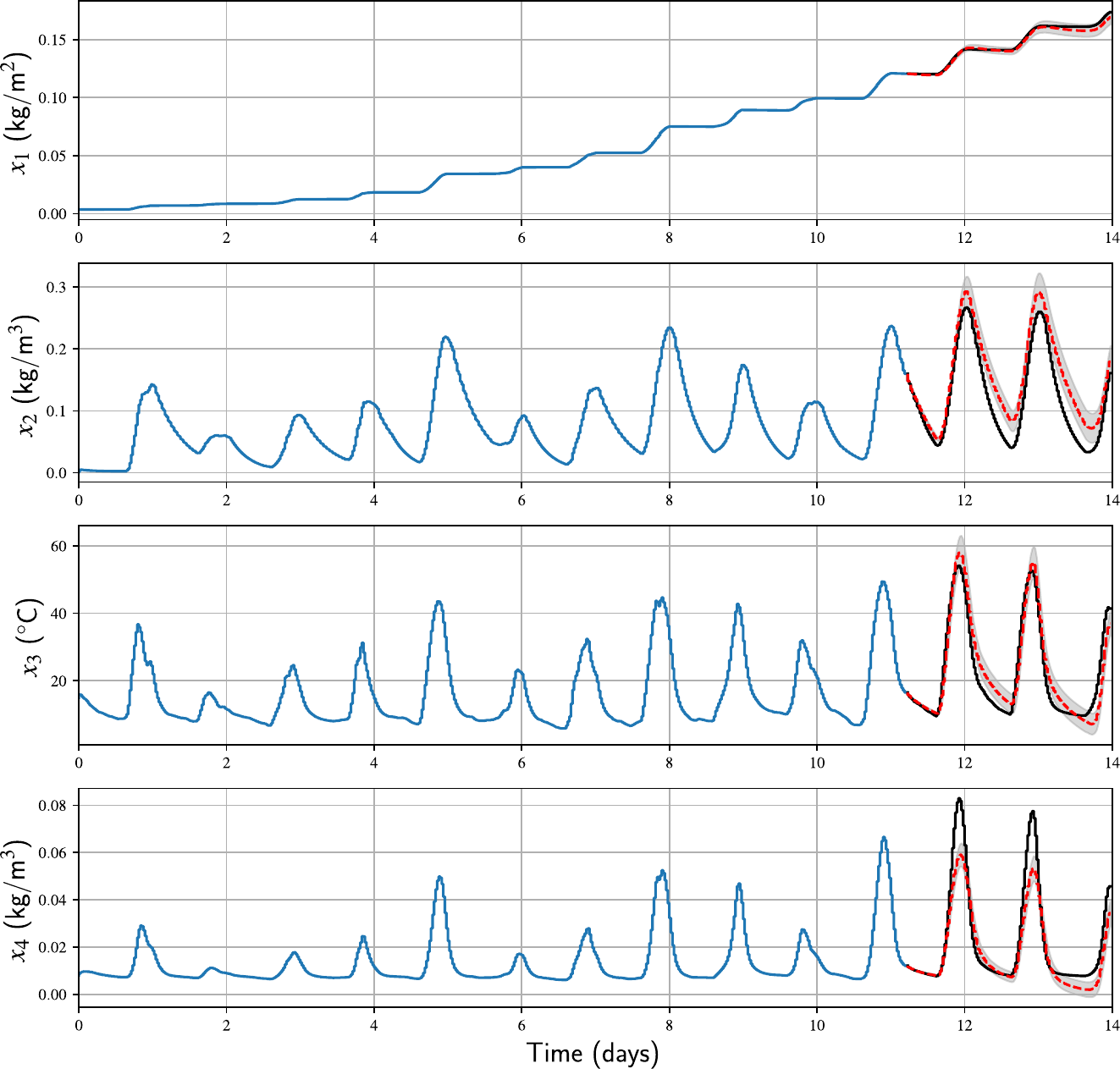}
	\caption{System identification results. The Bayesian neural ODE is trained using $N_s$ scenarios from day 0 to day 11. After day 11, the neural ODE is simulated $N_{\hat{p}}$ times for predicting 3 days ahead using~\eqref{eq:modelDiscrete1}. The red dashed curve is the mean of the predicted state variable and the gray area is its 99\% confidence interval. Black dashed is the true (unknown) state variable, which is shown for the comparison. \label{fig:sysIdResults}}
\end{figure}

Figure~\ref{fig:sysIdResults} shows the identification results for one scenario. The state variable for the first eleven days is generated with~\eqref{eq:system}. This specific time-series has not been used for training the Bayesian neural ODE. At around eleven days, the Bayesian neural ODE is propagated $N_{\hat{p}}$ times forward in time for three days. The initial conditions and weather prediction data are assumed to be known for these simulations. Subsequently, the mean state (red dashed in Fig.~\ref{fig:sysIdResults}) and 99\% confidence interval (gray area in Fig.~\ref{fig:sysIdResults}) are computed. 

For this specific time series, it can be concluded that the Bayesian neural ODE is capable of estimating the true (unknown) state variable (black dashed in Fig.~\ref{fig:sysIdResults}) for the upcoming three days. It should be noted however that initial conditions and weather are assumed to be known. The following step is therefore to investigate the effect of uncertain initial conditions and weather on the predictive capabilities of the Bayesian neural ODE. Also, models, as presented in this paper, can easily be used in stochastic predictive controllers as presented in~\cite{Boersma2022a}.

\clearpage
\section{CONCLUSIONS}
\label{sec:conclusions}

In this paper, we have presented a sparse Bayesian deep learning approach for identifying the dynamics of lettuce greenhouse models that incorporate both varying climate dynamics and crop growth status. The generated Bayesian neural ODE model is based on a deep neural network with the model parameters obeying certain probability distributions, rather than a single fixed value. The Bayesian inference approach has been applied to train the neural ODE using the data of multiple scenarios. Our simulation results have demonstrated that the developed model can capture the inherent nonlinear dynamics of the system, providing probabilistic prediction of environmental variables and lettuce growth within the greenhouse. The prediction result with uncertainty quantification provides a more comprehensive and realistic understanding of greenhouse dynamics, compared to the standard used deterministic models. This will also be potentially useful for more realistic control and decision-making in crop cultivation in horticulture.

\section*{ACKNOWLEDGMENT}
   The authors would like to thank Dr. Hongpeng Zhou from the University of Manchester for the insightful discussions.


\section*{APPENDIX}
The greenhouse with lettuce model is defined as:
\begin{equation*}
	\begin{aligned}
		\frac{\text{d} x_1(t)}{\text{d}t} &= p_{1} \phi_{\text{phot,c}}(t) -  p_{2} x_1(t) 2^{x_3(t)/10-5/2} , \\
		\frac{\text{d} x_2(t)}{\text{d}t} &= \frac{1}{p_{9}}   \Big( -\phi_{\text{phot,c}}(t) + p_{10} x_1(t)  2^{x_3(t)/10-5/2}  \\ 
		& \cdots + u_1(t) 10^{-6} - \phi_{\text{vent,c}}(t) \Big), \\
		\frac{\text{d} x_3(t)}{\text{d}t} &= \frac{1}{p_{16}} \Big( u_3(t) - (p_{17} u_2(t)10^{-3} +p_{18})  \\
		& \cdots (x_3(t)-d_3(t)) + p_{19} d_1(t) \Big), \\
		\frac{\text{d} x_4(t)}{\text{d}t} &= \frac{1}{p_{20}} \big( \phi_{\text{transp,h}}(t) - \phi_{\text{vent,h}}(t) \big) , \\
	\end{aligned}
\end{equation*}
with $t \in \mathbb{R}$ the continuous time and
\begin{equation*}
	\begin{aligned}
		\phi_{\text{phot,c}}(t) &= \Big(1-\text{exp} \big(-p_{3}x_1(t) \big) \Big) \\
		&\Big(p_{4}d_1(t)\big(-p_{5}x_3(t)^2+p_{6} x_3(t) - p_{7} \big) \\
		& \cdots \big( x_2(t)- p_{8}\big)\Big) /\varphi(t),\\
		\varphi(t)              &= p_{4}d_1(t)+\big(-p_{5}x_3(t)^2+p_{6} x_3(t) - \\
		& \cdots  p_{7} \big) \big( x_2(t)- p_{8}\big), \\
		\phi_{\text{vent,c}}(t) &= \big(u_2(t)10^{-3}+p_{11}\big)\big(x_2(t)-d_2(t)\big),\\
		\phi_{\text{transp,h}}(t) &= p_{21} \Big(1-\text{exp} \big(-p_{3}x_1(t) \big) \Big) \\ 
		\cdots  \Big( &\frac{p_{22}}{p_{23}(x_3(t)+p_{24})}  \text{exp} \Big( \frac{p_{25} x_3(t)}{x_3(t) + p_{26}} \Big) - x_4(t) \Big), \\
		\phi_{\text{vent,h}}(t) &= \big(u_2(t)10^{-3}+p_{11}\big)\big(x_4(t)-d_4(t)\big),		
	\end{aligned}
\end{equation*}
with $\phi_{\text{phot,c}}(t), \phi_{\text{vent,c}}(t),$ $\phi_{\text{transp,h}}(t)$ and $\phi_{\text{vent,h}}(t)$ are the gross canopy photosynthesis rate, mass exchange of CO$_2$ through the vents, canopy transpiration and mass exchange of H$_2$O through the vents, respectively. The initial state and control signal are: $x(0) = \begin{pmatrix} 0.0035 & 0.001 & 15 &0.008 \end{pmatrix}^T$, $u(0) = \begin{pmatrix} 0 & 0 & 0 \end{pmatrix}^T$. The model parameters $p_{i,j}$ are chosen following~\cite{vanHenten1994} and given in Table~\ref{tab:model_parameters}. The model is discretized using the explicit fourth order Runge-Kutta method resulting in the discrete-time model as presented in~\eqref{eq:model}. 

\begin{table}[h]
	\center
	\caption{\label{tab:model_parameters}Values of the model parameters that are taken from~\cite{vanHenten1994}.}
	\resizebox{.48\textwidth}{!}{	\begin{tabular}{ cc|cc|cc|cc } 
			parameter & value & parameter & value & parameter & value & parameter & value\\
			\hline
			$p_{1}$ & 0.544 & $p_{9}$ 	& 4.1 			& $p_{16}$ & 3$\cdot 10^{4}$ & $p_{20}$ 	& 4.1 \\ 
			$p_{2}$ & 2.65 $\cdot 10^{-7}$ & $p_{10}$ 	& 4.87 $\cdot 10^{-7}$ 			& $p_{17}$ & 1290 & $p_{21}$ & 0.0036 \\ 
			$p_{3}$ & 53 & $p_{11}$ 	& 7.5 $\cdot 10^{-6}$ 			& $p_{18}$ & 6.1 & $p_{22}$ 	& 9348 \\ 
			$p_{4}$ & 3.55 $\cdot 10^{-9}$ &	$p_{12}$	& 8.31  		& $p_{19}$ & 0.2 & $p_{23}$ 	& 8314 \\ 
			$p_{5}$ & 5.11 $\cdot 10^{-6}$ & $p_{13}$	& 273.15		&  			&   & $p_{24}$ 	& 273.15 \\ 
			$p_{6}$ & 2.3 $\cdot 10^{-4}$ & $p_{14}$	& 101325		&  			&   & $p_{25}$ 	& 17.4 \\ 
			$p_{7}$ & 6.29 $\cdot 10^{-4}$ & $p_{15}$	& 0.044			&  			&   & $p_{26}$ 	& 239 \\ 
			$p_{8}$ & 5.2 $\cdot 10^{-5}$ & 			&   			&  			&   &  $p_{27}$  	&   17.269  \\ 
			&   & 			&   			&  			&   							&  $p_{28}$	&  238.3  \\
	\end{tabular}}
\end{table}	

\bibliographystyle{ieeetr}
\bibliography{mybibfile}

\end{document}